\documentclass[12 pt,a4paper]{article} 

\usepackage[utf8]{inputenc}
\usepackage[english]{babel}  
\usepackage{graphicx}
\usepackage{lmodern}
\usepackage{amsmath}
\usepackage{tikz}
\usetikzlibrary{calc}
\usepackage{subfig}
\usepackage{longtable}
\usepackage{epsfig}
\usepackage{color}
\usepackage{setspace}
\usepackage[top=2cm,bottom=2.5cm,left=2.5cm,right=2.5cm]{geometry}
\usepackage{sectsty}
\usepackage{fancyhdr}
\usepackage{wasysym}
\usepackage[colorlinks=true]{hyperref}
\usepackage{caption}
\usetikzlibrary{patterns}
\usepackage{amsfonts}
\usepackage{amssymb}
\usepackage{amsthm}
\newtheorem{theorem}{Theorem}
\newtheorem{remark}[theorem]{Remark}
\newcommand{\sym}{\mathrm{Sym}}
\newcommand{\Skew}{\mathrm{Skew}}
\newcommand{\Div}{\mathrm{div}}
\newcommand{\tr}{\mathrm{tr}}
\newcommand{\one}{\mathbf{1}}

\title{A variational discrete element method for the computation of Cosserat elasticity}
\author{\begin{minipage}{\textwidth}\centering
		Frédéric Marazzato \\
		\small{Department of Mathematics, Louisiana State University, Baton Rouge, LA 70803, USA}\\
   \small{email: \texttt{marazzato@lsu.edu}\\}
   \end{minipage}
   }

\date{}   
\begin{document}

\maketitle

\begin{abstract}
The variational discrete element method developed in \cite{marazzato2020variational} for dynamic elasto-plastic computations is adapted to compute the deformation of elastic Cosserat materials.
In addition to cellwise displacement degrees of freedom (dofs), cellwise rotational dofs are added.
A reconstruction is devised to obtain $P^1$ non-conforming polynomials in each cell and thus constant strains and stresses in each cell.
The method requires only the usual macroscopic parameters of a Cosserat material and no microscopic parameter.
Numerical examples show the robustness of the method for both static and dynamic computations in two and three dimensions.
\end{abstract}

\section{Introduction}
Cosserat continua have been introduced in \cite{cosserat1909theorie}.
They generalize Cauchy continua by adding a miscroscopic rotation to every infinitesimal element. Cosserat continua can be considered as a generalization of Timoshenko beams to two and three-dimensional structures.
Contrary to traditional Cauchy continua of order one, Cosserat continua are able to reproduce some effects of the micro-structure of a material through the definition of a characteristic length written $\ell$ \cite{forest2001asymptotic}.
Cosserat media can appear as homogenization of masonry structures \cite{stefanou2008three,godio2017limit} or be used to model liquid crystals \cite{ericksen1974liquid}, Bingham--Cosserat fluids \cite{shelukhin2013cosserat} and localization in faults under shear deformation in rock mechanics \cite{RATTEZ201854}, for instance.

Discrete Element methods (DEM) have been introduced in \cite{hoover1974two} to model crystalline materials and in 
\cite{cundall1979discrete} for applications to geotechnical problems. Their use in granular materials and rock simulation is still widespread \cite{potyondy2004bonded,ries2007shear}.
Although DEM are able to represent accurately the behaviour of granular materials, their use to compute elastic materials is more delicate especially regarding the choice of microscopic material parameters.
The macroscopic parameters like Young modulus and Poisson ratio are typically recovered from numerical experiments using the set of microscopic parameters \cite{andre2012discrete,jebahi2015discrete}.
To remedy this problem couplings of DEM with Finite Element Methods (FEM) have been devised \cite{michael2015fem,avci2012fem}.
Beyond DEM-FEM couplings, attempts to simulate continuous materials with DEM have been proposed. In \cite{andre2012discrete,andre2013using}, the authors used a stress reconstruction inspired by statistical physics but the method suffers from the non-convergence of the macroscopic parameters with respect to the microscopic parameters.
In \cite{notsu2014symmetry} the authors derive a DEM method from a Lagrange $P^1$ FEM but cannot simulate materials with $\nu \geq 0.3$. In \cite{LM_CM_2012}, the authors pose the basis of variational DEM by deriving forces from potentials and link their method to Cosserat continua. Following this work, \cite{marazzato2020variational} proposed a variational DEM that can use polyhedral meshes and which is a full discretization of dynamic elasto-plasticity equations for a Cauchy continua.
However, in the process of approximating Cauchy continua, the unknown rotations in elements and torque between elements from \cite{LM_CM_2012} had to be removed.
The present paper builds on the achievements of \cite{marazzato2020variational} and reintroduces rotations by adding cellwise rotational degrees of freedom (dofs) to take into account micro-rotations thus leading to the natural discretization of Cosserat continua in lieu of Cauchy continua. Consequently it also greatly simplifies the integration of rotations in a dynamic evolution compared to \cite{LM_CM_2012} where a non-linear problem is solved at every time-step using an explicit RATTLE scheme.
Using the proposed method allows to use the usual tools of FEM applied to a DEM and thus makes its use and analysis much easier.
Also, the restrictions on meshes are relaxed allowing to use simplicial meshes in lieu of Voronoi meshes which are used in \cite{LM_CM_2012,andre2019novel} and are cumbersome to produce. 
Cosserat elasticity is usually computed through a $P^2$-$P^1$ Lagrange mixed element \cite{providas2002finite}.
However, the coupling with a traditional DEM is not obvious due to the location of the dofs in the methods (nodes for FEM and cell barycentres for DEM). In the proposed method, the dofs are located at the cell barycentres and thus a coupling with traditional DEM would be greatly simplified. Designing such a coupling is left for future work.
DEM can also be very useful in computing crack propagation due to their natural ability to represent discontinuous fields. In \cite{marazzato2021quasistatic}, the authors have built a variational DEM derived from \cite{marazzato2020variational} to compute cracking in elastic Cauchy materials. DEM are also very efficient in their ability to compute fragmentation \cite{puscas2015conservative}.
The proposed DEM can be considered as a first step towards its extension to compute cracking in elastic Cosserat materials.

In the present method, displacement and rotational dofs are placed at the barycentre of every cell and the Dirichlet boundary conditions are imposed weakly similarly to discontinuous Galerkin methods \cite{arnold1982interior}.
The general elasto-dynamic problem in a Cosserat continuum is presented in Section \ref{sec:continuous problem}. The discrete setting is then presented in Section \ref{sec:space discretization} and
alongside in Section \ref{sec:DEM interpretation}, the discrete strain-stress system derived from the continuous equations is reinterpreted in a DEM fashion as a force-displacement system.
In Section \ref{sec:numerical results}, validation tests are performed to prove the robustness and the precision of the chosen approach. In Section \ref{sec:dynamics}, the fully discrete system in space and time is described to compute dynamic evolutions and two and three dimensional test cases are presented proving the the correct integration of Cosserat elasticity.
Finally, Section \ref{sec:conclusion} draws some conclusions and presents potential subsequent work.

\section{Governing equations}
\label{sec:continuous problem}
An elastic material occupying, in the reference configuration, the domain $\Omega \subset \mathbb{R}^d$, where $d=2,3$, is considered to evolve dynamically over the time-interval $(0,T)$, where $T>0$, under the action of volumetric forces and boundary conditions.
The strain regime is limited to small strains and the material is supposed to have a micro-structure responding to a Cosserat material law. The material is also supposed to be isotropic and homogeneous.
Consequently, the material law is restricted to isotropic Cosserat elasticity in the following.
The displacement field is written $u \in \mathbb{R}^3$ and the rotation $R \in SO_3(\mathbb{R})$. Under the small strain regime, the rotation $R$ can be mapped to a micro-rotation $\varphi \in \mathbb{R}^3$ such that $\varphi = - \frac12 \epsilon : R$ and $R = \mathbf{1} - \epsilon \cdot \varphi$, where $\mathbf{1}$ is the unit three dimensional tensor and $\epsilon$ is a third-order tensor giving the signature of a permutation $(i,j,k)$. Thus $\epsilon_{ijk} = 1$ for an even permutation, $-1$ for an odd permutation and $0$, otherwise. For instance $\epsilon_{123} = 1$, $\epsilon_{112} = 0$ and $\epsilon_{132} = -1$.

\begin{remark} [2d case]
In two space dimensions, the micro-rotation is just a scalar $\varphi \in \mathbb{R}$ and $\epsilon$ is a matrix such that $\epsilon = \begin{pmatrix} 0 & 1 \\ -1 & 0 \\\end{pmatrix}$.
\end{remark}

The deformation tensor $e$ and the tension-curvature tensor $\kappa$ are defined as
\begin{equation}
\left\{
\begin{aligned}
& e(u,\varphi) = \nabla u + \epsilon \cdot \varphi \\
& \kappa(\varphi) = \nabla \varphi
\end{aligned}
\right.
\label{eq:cosserat continuous}
\end{equation}
The force-stress tensor $\sigma$ and the couple-stress tensor $\mu$ are linked to the strains $(e,\kappa)$ by
\begin{equation}
\left\{
\begin{aligned}
& \sigma(u,\varphi) = \mathbb{C} : e(u,\varphi) 	 \\
& \mu(\varphi) = \mathbb{D} : \kappa(\varphi),
\end{aligned}
\right.
\label{eq:loi comportement}
\end{equation}
where $\mathbb{C}$ and $\mathbb{D}$ are fourth-order tensors translating the material behaviour. Note that $\sigma(u,\varphi)$ is generally not symmetric unlike in Cauchy continua. 
In the three dimensional case $d=3$, we write:
\begin{equation}
\label{eq:material law}
\begin{aligned}
& \sigma(u,\varphi) = K \tr(e)(u,\varphi) \one + 2G \left( \sym(e)(u,\varphi) - \frac1{d} \tr(e)(u,\varphi) \one \right) + 2G_c \Skew(e)(u,\varphi), \\
& \mu(\varphi) = L \tr(\kappa)(\varphi) \one + 2M \left( \sym(\kappa)(\varphi) - \frac1{d} \tr(\kappa)(\varphi) \one \right)+ 2M_c \Skew(\kappa)(\varphi), \\
\end{aligned}
\end{equation}
where $K,G,G_c,L,M$ and $M_c$ are elastic moduli and $\sym$ gives the symmetric part of a rank-two tensor whereas $\Skew$ gives its skew-symmetric part. Under supplementary assumptions, the number of elastic moduli can be reduced from six to four in three dimensions \cite{jeong2010existence} but that is not the path chosen in this paper.
\begin{remark} [Length scale]
\label{rk:length}
Following \cite{forest2001asymptotic}, a characteristic length $\ell$ can be defined in a Cosserat elastic medium as $\max_{ijkl} \mathbb{C}_{ijkl} = \ell^2 \left(\max_{ijkl} \mathbb{D}_{ijkl} \right)$. This length can be interpreted as the length of the microstructure of the material. When $\ell \to 0$, the Cosserat material law can be homogenized into a Cauchy law \cite{forest2001asymptotic}.
\end{remark}
We introduce the volumic mass $\rho \in \mathbb{R}$ and the micro-inertia per unit mass $I \in \mathbb{R}$.
The dynamics equation in strong form write
\begin{equation}
\left\{
\begin{aligned}
& \Div(\sigma(u,\varphi)) + f = \rho \ddot{u}	 \\
& \Div(\mu(\varphi)) - \epsilon : \sigma(u,\varphi) + \mathfrak{c} = \rho I \ddot{\varphi}.
\end{aligned}
\right.
\label{eq:strong dynamics equation}
\end{equation}
Let $\partial \Omega = \partial \Omega_N \cup \partial \Omega_D$ be a partition of the boundary of $\Omega$. By convention $\partial \Omega_D$ is a closed set and $\partial \Omega_N$ is a relatively open set in $\partial \Omega$. The boundary $\partial \Omega_D$ has imposed displacements and micro-rotations $(u_D,\varphi_D)$, we thus enforce
\begin{equation}
\label{eq:Dirichlet boundary conditions}
\left\{
\begin{aligned}
& u = u_D \text{ on } \partial \Omega_D, \\
& \varphi = \varphi_D \text{ on } \partial \Omega_D.
\end{aligned}
\right.
\end{equation}
The normal and couple stresses $(g,m)$ are imposed on $\partial \Omega_N$, that is, we enforce
\begin{equation}
\label{eq:Neumann boundary conditions}
\left\{
\begin{aligned}
& \sigma \cdot n = g \text{ on } \partial \Omega_N, \\
& \mu \cdot n = m \text{ on } \partial \Omega_N.
\end{aligned}
\right.
\end{equation}
To write a variational DEM, we write the dynamics equations in weak form. Taking $(\tilde{u},\tilde{\varphi})$ as test functions, verifying homogeneous Dirichlet boundary conditions on $\partial \Omega_D$ ($\tilde{u}_{|\partial \Omega_D}=0 = \tilde{\varphi}_{|\partial \Omega_D}$), one has over $(0,T)$
\begin{multline}
\label{eq:weak form dynamics}
\int_{\Omega} \rho \ddot{u} \cdot v + \rho I \ddot{\varphi} \cdot \psi + \int_{\Omega} e(u,\varphi) : \mathbb{C} : e(v,\psi) + \kappa(\varphi) : \mathbb{D} : \kappa(\psi) \\ =  \int_{\Omega} f \cdot v + \mathfrak{c} \cdot \psi + \int_{\partial \Omega_N} g \cdot v + m \cdot \psi,
\end{multline}
while still imposing the Dirichlet boundary conditions of Equation \eqref{eq:Dirichlet boundary conditions}. Note that the bilinear form in the left-hand side of \eqref{eq:weak form dynamics} is symmetric.
\begin{remark}
Initial conditions on $(u,\varphi)$ and $(\dot{u},\dot{\varphi})$ need to be specified to compute the solution of Equation \ref{eq:weak form dynamics}.
\end{remark}

\section{Space semidiscretization}
\label{sec:space discretization}
The domain $\Omega$ is discretized with a
mesh $\mathcal{T}_h$ of size $h$ made of polyhedra with planar facets in three space dimensions or polygons with straight edges in two space dimensions. 
We assume that $\Omega$ is itself a polyhedron or a polygon so that the mesh covers $\Omega$ exactly, and we also assume that the mesh is compatible with the partition of the boundary $\partial\Omega$ into the Dirichlet and Neumann parts. 

\subsection{Degrees of freedom}
\begin{figure} [!htp] 
\begin{center}
\begin{tikzpicture}[scale=1.3]
\coordinate (a) at (0,-1);
\coordinate (b) at (0,1);
\coordinate (c) at (1.3,1);
\coordinate (d) at (-1.7,-0.8);
\coordinate (e) at (-2,-3);
\coordinate (f) at (0.5,-2.5);
\coordinate (g) at (1.7,0);
\coordinate (h) at (2.8,0.1);
\coordinate (i) at (2.1,-1.1);
\coordinate (j) at (3.5,-1.5);
\coordinate (k) at (2.2,-3.5);
\coordinate (l) at (0.8,-4);
\coordinate (m) at (-1.7,-4);

\coordinate (abd) at (barycentric cs:a=0.3,b=0.3,d=0.3);

\coordinate (ade) at (barycentric cs:a=0.3,e=0.3,d=0.3);
\coordinate (aef) at (barycentric cs:f=0.3,e=0.3,a=0.3);
\coordinate (abcg) at (barycentric cs:a=0.25,g=0.25,c=0.25,b=0.25);
\coordinate (cgh) at (barycentric cs:g=0.3,h=0.3,c=0.3);
\coordinate (afgi) at (barycentric cs:a=0.25,g=0.25,i=0.25,f=0.25);
\coordinate (ghij) at (barycentric cs:h=0.25,g=0.25,i=0.25,j=0.25);
\coordinate (fijkl) at (barycentric cs:i=0.2,j=0.2,k=0.2,l=0.2,f=0.2);
\coordinate (eflm) at (barycentric cs:f=0.25,l=0.25,m=0.25,e=0.25);

\draw[fill=gray,opacity=0.1] (b) --(c) -- (h) -- (j) -- (k) -- (l) --
(m) -- (e) -- (d) -- cycle;
\node[right] at (barycentric cs:f=0.8,a=0.25,g=0.25,i=0.25) {\Large
  $\Omega$};
\node[right] at (barycentric cs:k=0.5,j=0.5) {\Large $\partial\Omega$};
\path[draw] (a)-- (b)-- (c) -- (g) -- cycle;
\path[draw] (a)-- (b)-- (d) -- cycle;
\path[draw] (a) -- (e) -- (d) -- cycle;
\path[draw] (a) -- (f) -- (e) -- cycle;
\path[draw] (a) -- (g) -- (i) -- (f) -- cycle;
\path[draw] (c) -- (h) -- (g) -- cycle;
\path[draw] (g) -- (h) -- (j) -- (i) -- cycle;
\path[draw] (i) -- (j) -- (k) -- (l) -- (f) -- cycle;
\path[draw] (f) -- (l) -- (m) -- (e) -- cycle;

\node[blue,right] at (abd) {$\blacktriangledown$};
\node[blue,right] at (ade) {$\blacktriangledown$};
\node[blue,right] at (afgi) {$\blacktriangledown$};
\node[blue,right] at (ghij) {$\blacktriangledown$};
\node[blue,right] at (fijkl) {$\blacktriangledown$};
\node[blue,right] at (eflm) {$\blacktriangledown$};
\node[blue,right] at (aef) {$\blacktriangledown$};
\node[blue,right] at (cgh) {$\blacktriangledown$};
\node[blue,right] at (abcg) {$\blacktriangledown$};

\node[red] at (abd) {$\blacktriangle$};
\node[red] at (abd) {$\blacktriangle$};
\node[red] at (ade) {$\blacktriangle$};
\node[red] at (afgi) {$\blacktriangle$};
\node[red] at (ghij) {$\blacktriangle$};
\node[red] at (fijkl) {$\blacktriangle$};
\node[red] at (eflm) {$\blacktriangle$};
\node[red] at (aef) {$\blacktriangle$};
\node[red] at (cgh) {$\blacktriangle$};
\node[red] at (abcg) {$\blacktriangle$};

\node[red] at (4.6,-1) {$\blacktriangle$};
\node[right] at (4.7,-1) {$(u_c)_{c\in\mathcal{C}}$};
\node[blue] at (4.6,-2.05) {$\blacktriangledown$};
\node[right] at (4.7,-2) {$(\phi_c)_{c\in\mathcal{C}}$};
\end{tikzpicture}
\caption{Continuum $\Omega$ covered by a polyhedral mesh and vector-valued degrees of freedom for the displacement.}
\label{fig:dofs}
\end{center}
\end{figure}
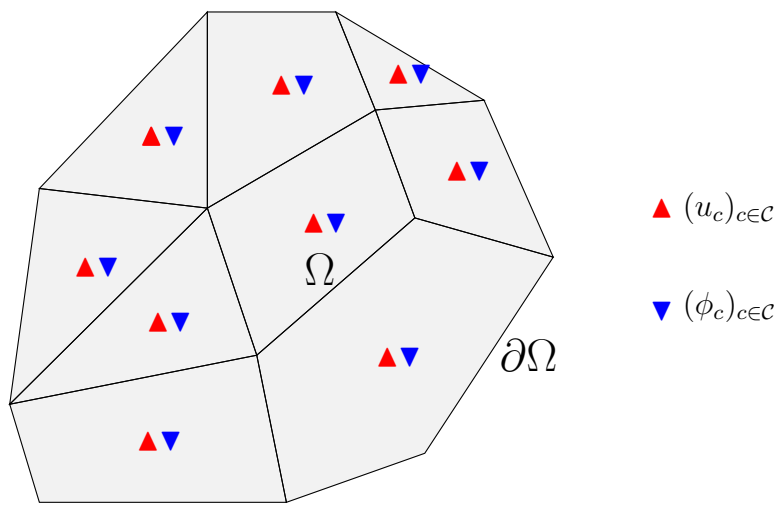

Let $\mathcal{C}$ denote the set of mesh cells. Pairs of vector-valued volumetric degrees of freedom (dofs) for a generic displacement field and a generic micro-rotation field $(v_h,\psi_h):=(v_c,\psi_c)_{c\in\mathcal{C}}\in\mathbb{R}^{2d\#(\mathcal{C})}$ are placed at the barycentre of every mesh cell $c\in\mathcal{C}$, where $\#(S)$ denotes the cardinality of any set $S$. Figure \ref{fig:dofs} illustrates the position of the dofs in the mesh.

Let $\mathcal{F}$ denote the set of mesh facets. We partition 
this set as $\mathcal{F}=\mathcal{F}^i\cup \mathcal{F}^b$, where 
$\mathcal{F}^i$ is the collection of the internal facets shared
by two mesh cells and $\mathcal{F}^b$ is the collection of the boundary facets 
sitting on the boundary $\partial\Omega$ (such facets belong to the boundary of only
one mesh cell). The sets $\mathcal{F}^b_N$ and $\mathcal{F}^b_D$ are defined as a partition of $\mathcal{F}^b$ such that $\forall F \in \mathcal{F}^b_N, F \subset \partial \Omega_N$ and $\forall F \in \mathcal{F}^b_D, F \subset \partial \Omega_D$.

\subsection{Facet reconstructions}
Using the cell dofs introduced above, we reconstruct a collection of
displacements and micro-rotations $(v_{\mathcal{F}},\psi_{\mathcal{F}}) :=(v_F,\psi_F)_{F\in\mathcal{F}}\in
\mathbb{R}^{2d\#(\mathcal{F})}$ on the mesh facets. 
The facet reconstruction operator is
denoted $\mathcal{R}$ and we write
\begin{equation}
(v_{\mathcal{F}},\psi_{\mathcal{F}}) := \left(\mathcal{R}(v_h),\mathcal{R}(\psi_h)\right).
\end{equation}

The reconstruction operator $\mathcal{R}$ is constructed in the same way as in the finite volume methods studied in \cite[Sec. 2.2]{eymard2009discretization} and the variational DEM developed in \cite{marazzato2020variational}.
For a given facet $F\in\mathcal{F}$, we select  neighbouring cells collected in a subset denoted $\mathcal{C}_F$, as well as coefficients $(\alpha_F^c)_{c\in\mathcal{C}_F}$ and we set
\begin{equation}\label{eq:barycentric}
\left\{
\begin{aligned}
&\mathcal{R}_F(v_h) := \sum_{c\in\mathcal{C}_F}\alpha_F^c v_c, \qquad \forall v_h\in V_h, \\
&\mathcal{R}_F(\psi_h) := \sum_{c\in\mathcal{C}_F}\alpha_F^c \psi_c, \qquad \forall \varphi_h\in V_h.
\end{aligned}
\right.
\end{equation}
The reconstruction is based on barycentric coordinates.
The coefficients $\alpha_F^c$ are chosen as the barycentric coordinates of the facet barycentre $\mathbf{x}_F$ in terms of the location of the barycenters of the cells in $\mathcal{C}_F$. For any facet $F \in \mathcal{F}$, the set $\mathcal{C}_F$ is constructed so as to contain exactly $(d+1)$ points forming the vertices of a non-degenerate simplex. Thus, the barycentric coefficients $(\alpha^c_F)_{c \in \mathcal{C}_F}$ are computed by solving the linear system:
\begin{equation}
\label{eq:system barycentric coordinates}
\left\{
\begin{alignedat}{2}
&\sum_{c\in\mathcal{C}_F}{\alpha_F^c} = 1,\qquad &\forall F\in\mathcal{F}, \\
&\sum_{c\in\mathcal{C}_F}{\alpha_F^c \mathbf{x}_c} = \mathbf{x}_F,\qquad &\forall F\in\mathcal{F}, \\
\end{alignedat}
\right.
\end{equation}
where $\mathbf{x}_c$ is the position of the barycenter of the cell $c$.
An algorithm is presented thereafter to explain the selection of the neighbouring dofs in $\mathcal{C}_F$.
This algorithm has to be viewed more as a proof-of-concept than as an optimized algorithm. For a more involved algorithm, see \cite{marazzato2020variational}.
We observe that this algorithm is only used in a preprocessing stage of the computations. For a given facet $F \in \mathcal{F}$,
\begin{enumerate}
\item Assemble in a set $\mathcal{N}_F$ the cell or cells containing the facet $F$. Then, add to $\mathcal{N}_F$ the cells sharing a facet with the cells already in $\mathcal{N}_F$. Repeat the last operation once.

\item Select a subset $\mathcal{C}_F$ of $\mathcal{N}_F$ with exactly $(d+1)$ elements and whose cell barycenters form a non-degenerate simplex (tetrahedron in 3d and triangle in 2d).
\end{enumerate}
This algorithm ensures that the dofs selected for the reconstruction in Equation \eqref{eq:barycentric} remain $\mathcal{O}(h)$ close to the facet $F$.
\begin{remark}
The last operation of step 1 described above which consists in enlarging the set $\mathcal{N}_F$ is generally not necessary in two space dimensions but it becomes necessary in three space dimensions as some locally complicated mesh geometries can cause all barycenters of cells in $\mathcal{N}_F$ to form a degenerate simplex. In that case, step 2 cannot be performed correctly.
\end{remark}

\begin{remark} [Influence of the choice of $\mathcal{C}_F$]
The choice of the elements in $\mathcal{C}_F$ for $F \in \mathcal{F}$ has an impact on the eigenvalues of the rigidity matrix and thus on its conditioning. This impact is explored in \cite{marazzato2020variational} regarding the CFL condition in explicit dynamic computations.
\end{remark}

\subsection{Gradient reconstruction}
Using the reconstructed facet displacements and micro-rotations and a discrete Stokes formula, it is possible to devise a discrete $\mathbb{R}^{d\times d}$-valued
piecewise-constant gradient field for the displacement and the micro-rotation
that we write $\mathcal{G}_{\mathcal{C}}(v_{\mathcal{F}}) := (\mathcal{G}_{c}(v_{\mathcal{F}}))_{c\in\mathcal{C}}
\in \mathbb{R}^{d^2\#(\mathcal{C})}$ and $\mathcal{G}_{\mathcal{C}}(\psi_{\mathcal{F}}) := (\mathcal{G}_{c}(\psi_{\mathcal{F}}))_{c\in\mathcal{C}}
\in \mathbb{R}^{d^2\#(\mathcal{C})}$.
Specifically we set in every mesh cell $c\in\mathcal{C}$,
\begin{equation}
\label{eq:gradient reconstruction}  
\left\{
\begin{aligned}  
&\mathcal{G}_c(v_{\mathcal{F}}) := \sum_{F \in \partial c} \frac{|F|}{|c|} v_F \otimes n_{F,c},
\qquad \forall v_{\mathcal{F}} \in\mathbb{R}^{d\#(\mathcal{F})}, \\
& \mathcal{G}_c(\psi_{\mathcal{F}}) := \sum_{F \in \partial c} \frac{|F|}{|c|} \psi_F \otimes n_{F,c},
\qquad \forall \psi_{\mathcal{F}} \in\mathbb{R}^{d\#(\mathcal{F})},
\end{aligned}
\right.
\end{equation}
where the summation is over the facets $F$ of $c$ and $n_{F,c}$ is the outward 
normal to $c$ on $F$. 
Consequently, the strains are defined for $c \in \mathcal{C}$ as
\begin{equation}
\left\{
\begin{aligned}
&e_{c}(v_h) := \mathcal{G}_{c}(v_h) + \epsilon \cdot \psi_c \in \mathbb{R}^{d\times d}, \\
& \kappa_c(\psi_h) := \mathcal{G}_c(\psi_h) \in \mathbb{R}^{d\times d},
\end{aligned}
\right.
\end{equation}
where $\mathcal{G}_c(v_h) := \mathcal{G}_c(\mathcal{R}(v_h))$ and $\mathcal{G}_c(\psi_h) := \mathcal{G}_c(\mathcal{R}(\psi_h))$.
Consequently, the discrete bilinear form of elastic energies writes
\begin{equation}
\label{eq:elastic bilinear form}
a_{\text{elas}}\left((u_h,\varphi_h);(v_h,\psi_h)\right) := \int_{\Omega} \left( e_h\left((u_h,\varphi_h)\right) : \mathbb{C} : e_h\left((v_h,\psi_h)\right) + \kappa_h(\varphi_h) : \mathbb{D} : \kappa_h(\psi_h) \right).
\end{equation}
Finally, we define two additional reconstructions which will be used to impose the Dirichlet boundary conditions. The reconstructions are written $\mathfrak{R}$ and consist in a collection of cellwise nonconforming $P^1$ polynomials defined for all $c \in \mathcal{C}$ by
\begin{equation}
\label{eq:DG reconstruction}
\left\{
\begin{aligned}
&\mathfrak{R}_c(v_h)(\mathbf{x}) := v_c + \mathcal{G}_c(v_h) \cdot (\mathbf{x} - \mathbf{x}_c), \\
& \mathfrak{R}_c(\psi_h)(\mathbf{x}) := v_c + \mathcal{G}_c(\psi_h) \cdot (\mathbf{x} - \mathbf{x}_c).
\end{aligned}
\right.
\end{equation}
where $\mathbf{x} \in c$ and $\mathbf{x}_c$ is the barycentre of the cell $c$.

\subsection{Mass bilinear form}
In the DEM spirit, the reconstruction of functions is chosen as constant in each cell so as to obtain a diagonal mass matrix.
The mass bilinear form is thus defined as
\begin{equation}
\label{eq:mass bilinear form}
m_h\left((u_h,\varphi_h);(v_h,\psi_h)\right) := \sum_{c \in \mathcal{C}} \rho |c| \left(u_c \cdot v_c + I\varphi_c \cdot \psi_c \right).
\end{equation}

\subsection{Discrete problem}
The discrete problem is defined as a lowest-order discontinuous Galerkin method similar to \cite{eymard2009discretization} and \cite{di2012cell}.
Consequently, penalty terms will be added to the discrete formulation of Equation \eqref{eq:cosserat continuous} for two reasons. The first is that the gradient reconstruction of Equation \eqref{eq:gradient reconstruction} cannot by itself control $(u_h,\phi_h)$ and thus a least-square penalty term will be added to the formulation to ensure the well-posedness of the problem.
More details can be found in \cite{eymard2009discretization} and \cite{di2012cell}.
The second is that the Dirichlet boundary conditions will not be imposed strongly but weakly through a non-symmetric Nitsche penalty acting on the boundary facets in $\partial \Omega_D$.

\subsubsection{Least-square penalty}

For an interior facet $F\in \mathcal{F}^i$, writing $c_{F,-}$ and $c_{F,+}$ the two mesh cells sharing $F$, that is, $F = \partial c_{F,-} \cap \partial c_{F,+}$, and orienting $F$ by the unit normal vector $n_F$ pointing from $c_{F,-}$ to $c_{F,+}$, we define
\begin{equation}
[\mathfrak{R}(v_h)]_F := \mathfrak{R}_{c_{F,-}}(v_h) - \mathfrak{R}_{c_{F,+}}(v_h).
\end{equation}
$[\mathfrak{R}(\psi_h)]_F$ is defined similarly.
The interior penalty term is defined as
\begin{multline}
\label{eq:innner penalty}
a_{\text{inner\_pen}}\left((u_h,\varphi_h);(v_h,\psi_h)\right) : = \sum_{F \in \mathcal{F}^i} \frac{1}{h_F} \int_F    ([\mathfrak{R}(u_h)]_F \otimes n_F) : \mathbb{C} : ( [\mathfrak{R}(v_h)]_F \otimes n_F) \\
+ ([\mathfrak{R}(\varphi_h)]_F \otimes n_F) : \mathbb{D} : ([\mathfrak{R}(\psi_h)]_F \otimes n_F).
\end{multline}

\subsubsection{Non-symmetric Nitsche penalty}
As the Dirichlet boundary conditions are imposed weakly, the discrete test functions do not verify $\mathcal{R}(v_h) = 0$ and $\mathcal{R}(\psi_h) = 0$ on $\partial \Omega_D$. Thus the following term, called the consistency term, coming from the integration by parts, leading from \eqref{eq:strong dynamics equation} to \eqref{eq:weak form dynamics}, must be taken into account
\begin{equation}
\label{eq:nitsche con}
a_{\text{con}}\left((u_h,\varphi_h);(v_h,\psi_h)\right) := - \int_{\partial \Omega_D} (\sigma_h(u_h,\varphi_h) \cdot n) \cdot \mathcal{R}(v_h) + (\mu_h(\varphi_h) \cdot n) \cdot \mathcal{R}(\psi_h). 
\end{equation}
Following ideas from \cite{burman2012penalty,boiveau2016penalty}, we introduce the following non-symmetric term, to obtain a method that is stable without having to add a least-square penalty term on $\partial \Omega_D$,
\begin{equation}
a_{\text{nsym}}\left((u_h,\varphi_h);(v_h,\psi_h)\right) := \int_{\partial \Omega_D} (\sigma_h(v_h,\psi_h) \cdot n) \cdot \mathcal{R}(u_h) + (\mu_h(\psi_h) \cdot n) \cdot \mathcal{R}(\varphi_h). 
\label{eq:nitsche nsym}
\end{equation}
A corresponding linear form is added to compensate the term in \eqref{eq:nitsche nsym} when the Dirichlet boundary conditions are verified exactly,
\begin{equation}
l_{\text{nsym}}\left((u_h,\varphi_h);(v_h,\psi_h)\right) := \int_{\partial \Omega_D} (\sigma_h(v_h,\psi_h) \cdot n) \cdot u_D + (\mu_h(\psi_h) \cdot n) \cdot \varphi_D.  
\label{eq:nitsche lin}
\end{equation}

\subsubsection{Discrete problem}
The bilinear form $a_h$ is defined as $a_h := a_{\text{elas}} + a_{\text{inner\_pen}} + a_{\text{con}} + a_{\text{nsym}}$.
The discrete problem may then be written: search for $(u_h,\varphi_h)$ such that for all $(v_h,\psi_h)$, one has
\begin{equation}
\label{eq:discrete dynamics}
m_h\left((\ddot{u}_h,\ddot{\varphi}_h);(v_h,\psi_h)\right) + a_h\left((u_h,\varphi_h);(v_h,\psi_h)\right) = l_h(v_h,\psi_h) + l_D(v_h,\psi_h) + l_\text{nsym}(v_h,\psi_h),
\end{equation}
where $l_h$ is the linear form that takes into account Neumann boundary conditions and volumic loads. It is defined as
\begin{equation}
\label{eq:forces exterior loads}
l_h(v_h,\psi_h) := \sum_{c \in \mathcal{C}} \left(\int_c f \right) \cdot v_c + \left( \int_c \mathfrak{c} \right) \cdot \psi_c + \sum_{F \in \mathcal{F}^b_N} \left(\int_F g \right) \cdot v_{c_F} + \left( \int_F m \right) \cdot \psi_{c_F}. 
\end{equation}

\subsection{Interpretation as a DEM}
\label{sec:DEM interpretation}
Traditional DEM are force-displacement systems in the sense that the deformation of the domain $\Omega$ is computed through the displacement of particles interacting through nearest-neighbours forces \cite{jebahi2015discrete}.
A major difference appears with the proposed DEM in which the interactions are not nearest-neighbours but have a larger stencil due to the facet and gradient reconstructions \eqref{eq:barycentric} and \eqref{eq:gradient reconstruction}.
\begin{figure} [!htp]
\centering
\subfloat{
\begin{tikzpicture}
\coordinate (a) at (0,0);
\coordinate (b) at (1,0);
\coordinate (c) at (0,1);
\coordinate (d) at (0,-1);
\coordinate (e) at (-1,0);
\fill[red] (a) circle (2mm);
\fill[blue] (b) circle (2mm);
\fill[blue] (c) circle (2mm);
\fill[blue] (d) circle (2mm);
\fill[blue] (e) circle (2mm);
\coordinate (f) at (1,1);
\coordinate (g) at (1,-1);
\coordinate (h) at (-1,-1);
\coordinate (i) at (-1,1);
\fill[green] (f) circle (2mm);
\fill[green] (g) circle (2mm);
\fill[green] (h) circle (2mm);
\fill[green] (i) circle (2mm);

\fill[red] (3,1) circle (2mm);
\node at (5,1) {Reference particle};
\fill[blue] (3,0) circle (2mm);
\node at (5.3,0) {Neighbouring particle};
\fill[green] (3,-1) circle (2mm);
\node at (5.5,-1) {Neighbour to Neighbour};
\end{tikzpicture}
}
\subfloat{
\begin{tikzpicture}
\coordinate (a) at (0,0);
\coordinate (b) at (1,0);
\coordinate (c) at (0,1);
\coordinate (abc) at (barycentric cs:a=0.33,b=0.33,c=0.33);
\fill[red] (abc) circle (2mm);
\draw (a) -- (b) -- (c) -- cycle;

\coordinate (d) at (0,-1);
\draw (b) -- (d) -- (a);
\coordinate (abd) at (barycentric cs:a=0.33,b=0.33,d=0.33);
\draw[blue, fill=blue] (abd) circle (2mm);

\coordinate (e) at (-1,0);
\draw (c) -- (e) -- (a);
\coordinate (ace) at (barycentric cs:a=0.33,c=0.33,e=0.33);
\fill[blue] (ace) circle (2mm);

\coordinate (f) at (1,1);
\draw (b) -- (f) -- (c);
\coordinate (bcf) at (barycentric cs:b=0.33,c=0.33,f=0.33);
\fill[blue] (bcf) circle (2mm);

\coordinate (g) at (0,2);
\draw (f) -- (g) -- (c);
\coordinate (cfg) at (barycentric cs:c=0.33,f=0.33,g=0.33);
\draw[green, fill=green] (cfg) circle (2mm);

\draw (d) -- (e);
\coordinate (ade) at (barycentric cs:a=0.33,d=0.33,e=0.33);
\fill[green] (ade) circle (2mm);

\coordinate (h) at (1,-1);
\draw (d) -- (h) -- (b);
\coordinate (bdh) at (barycentric cs:b=0.33,d=0.33,h=0.33);
\fill[green] (bdh) circle (2mm);

\draw[dashed] (abc) -- (bcf) -- (cfg) -- cycle;
\draw[dashed] (abc) -- (abd) -- (bdh) -- cycle;
\draw[dashed] (abc) -- (ace) -- (ade) -- cycle;

\draw[dashed] (2,0) -- (2.5,0);
\node at (3.5,0.3) {For facet};
\node at (4,-0.3) {reconstruction};
\end{tikzpicture}
}
\caption{DEM interpretation: left: traditional DEM. right: proposed variational DEM.}
\label{fig:DEM interactions}
\end{figure}
As can be seen in Figure \ref{fig:DEM interactions} on the left, in a traditional DEM, a particle (represented in red) interacts with its closest-neighbours (in blue) but not the neighbours of its neighbours (in green).
On the contrary, such interactions are present with the the proposed method as can be seen in Figure \ref{fig:DEM interactions} on the right.
In traditional DEM, the force and torque between two particles can be parametrized, for instance, by elastic springs in tension, shear and torque, or by beam elements \cite{jebahi2015discrete}.
The parameters of these springs or beams are called microscopic parameters.
Unfortunately, calibrating the microscopic parameters to simulate a given continuous material is difficult \cite{jeong2010existence}.
As the elastic bilinear form \eqref{eq:elastic bilinear form} is written in terms of strains and stresses, as in continuous materials, forces and torques are not explicit.
However, it is possible to rewrite \eqref{eq:elastic bilinear form} so as to extract forces between neighbouring cells considered as discrete elements.
The major advantage being that the forces retrieved are parametrized only by the continuous material parameters and not by microscopic parameters. Let us do so in the following.

The average value in an inner facet $F \in \mathcal{F}^i$ of a quantity $a$ is defined as $\{a\}_F := \frac12 (a_{c_{F,-}} + a_{c_{F,+}})$.
Rewriting Equation \eqref{eq:elastic bilinear form} and neglecting second order terms, one has
\begin{equation}
\label{eq:forces DEM}
\begin{alignedat}{2}
-a_{\text{elas}}((u_h,\varphi_h);(v_h,\psi_h)) & = & \sum_{F \in \mathcal{F}^i} |F| (\{\sigma_h(u_h, \varphi_h) \}_F \cdot n_F) \cdot [v_h]_F \\
& + & \sum_{F \in \mathcal{F}^i} |F| (\{\mu_h(\varphi_h) \}_F \cdot n_F) \cdot [\psi_h]_F \\
& + & \sum_{c \in \mathcal{C}} |c| (\epsilon : \sigma_c(u_h, \varphi_h)) \cdot \psi_c \\
& + & \sum_{F \in \mathcal{F}^b} |F| (\sigma_{c_F}(u_h, \varphi_h) \cdot n_F) \cdot (v_{c_F} - \mathcal{R}_F(v_h)) \\
& + & \sum_{F \in \mathcal{F}^b} |F| (\mu_{c_F}(\varphi_h) \cdot n_F) \cdot (\psi_{c_F} - \mathcal{R}_F(\psi_h)) \\
& + & \mathcal{O}(h^2),
\end{alignedat}
\end{equation}
where $c_F$ designates the unique cell containing a boundary facet $F \in \mathcal{F}^b$, $\sigma_h(u_h,\varphi_c) := \mathbb{C} : e_h(u_h,\varphi_c)$, $\mu_h(\varphi_h) := \mathbb{D} : \kappa_h(\varphi_h)$ and for an inner facet $F$, $[v_h]_F := v_{c_{F,-}} - v_{c_{F,+}}$ and $[\psi_h]_F := \psi_{c_{F,-}} - \psi_{c_{F,+}}$. 
The principle of action-reaction (or Newton's third law) can be read in the first two lines of the equation through the action of jump terms.
The third line represents the work of the momentum coming from the stresses. The fourth and fifth lines represent the work of the internal forces and momenta but related to boundary facets.
One can write the dynamics equations of an interior cell $c \in \mathcal{C}$ (or discrete element), with no facet on the boundary, as follows
\begin{equation}
\left\{
\begin{aligned}
& \rho |c| \ddot{u}_c \simeq \sum_{F \in \mathcal{F}^i, F \subset \partial c} \iota_{c,F} |F| \{\sigma_h(u_h, \varphi_h) \}_F \cdot n_F + \int_c f, \\ 
& \rho |c| I \cdot \ddot{\varphi}_c \simeq |c| \epsilon : \sigma_c(u_h, \varphi_h) + \sum_{F \in \mathcal{F}^i, F \subset \partial c} \iota_{c,F} |F| \{\mu_h(\varphi_h) \}_F \cdot n_F + \int_c \mathfrak{c}, \\
\end{aligned}
\right.
\end{equation}
up to second order terms and penalty terms and where $\iota_{c,F} = 1$ if $c = c_{F,-}$ and $\iota_{c,F} = -1$ if $c = c_{F,+}$. Each facet thus represents a link between two discrete elements and the forces and momenta are average quantities computed from cell-wise stress and momentum reconstructions. To obtain similar equations for a cell having a facet on $\partial \Omega$, one can refer to \cite{marazzato2020variational}.

\subsection{Implementation}
The method has been implemented in Python and is available at \url{https://github.com/marazzaf/DEM_cosserat.git.}
A finite element implementation available in \cite{sautot2014extension} has been used as a foundation for the implementation of the previous method.
FEniCS \cite{logg2012automated,10.1145/1731022.1731030} is used to handle meshes, compute facet reconstructions and assemble matrices.
PETSc \cite{DALCIN20111124,balay2019petsc} is used for matrix storage, matrix operations and as a linear solver.
Because FEniCS only supports simplicial meshes, the meshes used with the implementation are only triangular or tetrahedral. However, the method supports general polyhedra.

\subsection{Validation test cases}
\label{sec:validation test}
The validation test cases are two-dimensional. Cosserat elasticity is characterized only by four parameters in two space dimensions.
Following \cite{providas2002finite}, we chose the material parameters to be $G$ the shear modulus, $\ell$ the characteristic length of the microstructure, $a$ (which measures the ratio $G_c$ over $G$) and $\nu$ the Poisson ratio.
Let $\mathfrak{a} = \frac{2(1-\nu)}{1-2\nu}$ and $\mathfrak{b} = \frac{2\nu}{1-2\nu}$.
The material law then writes
\begin{equation}
\label{eq:2d material law}
\left( \begin{array}{c} \sigma_{xx} \\ \sigma_{yy} \\ \sigma_{xy} \\ \sigma_{yx} \end{array} \right)
= G \begin{pmatrix}
& \mathfrak{a} & \mathfrak{b} & 0 & 0 \\
& \mathfrak{b} & \mathfrak{a} & 0 & 0 \\
& 0 & 0 & 1+a & 1-a \\
& 0 & 0 & 1-a & 1+a \\
\end{pmatrix}
\cdot
\left( \begin{array}{c}
e_{xx} \\ e_{yy} \\ e_{xy} \\ e_{yx} \\
\end{array} \right),
\end{equation}
and
\begin{equation}
\label{eq:2d material law 2}
\mu = 4G\ell^2 \kappa.
\end{equation}
The following test has been used in \cite{providas2002finite} to validate the P2-P1 finite element in two space dimension. It consists in three computations on the rectangular domain $[-0.12,0.12]\times[0,0.12]$ with material parameters $G=10^3$Pa, $\ell=0.1$m, $a=0.5$ and $\nu=0.25$.
Estimates of the condition number of the rigidity matrices with and without the terms corresponding to equation \eqref{eq:innner penalty} are also computed.
The condition number of the rigidity matrix is approximated using \cite{fong2011lsmr} implemented in \texttt{scipy.sparse.linalg} which is a module of Scipy \cite{2020SciPy-NMeth}.
A single structured mesh containing $2,500$ elements is used to compute the condition numbers for the three test cases. It corresponds to $7,500$ dofs for the proposed DEM.

\subsubsection{First patch test}
The solution to the first test is 
\begin{equation}
\left\{
\begin{aligned}
& u_x(x,y) = \frac1{G}(x + \frac{y}2), \\
& u_y(x,y) = \frac1{G}(x + y), \\
&\varphi(x,y) = \frac1{4G}.
\end{aligned}
\right.
\end{equation}
Dirichlet boundary conditions are imposed on the entire boundary.
The volumic load is $f_x = f_y = \mathfrak{c} = 0$.
The condition number for the DEM with the penalty term \eqref{eq:innner penalty} is $83$ and $51$ without it.
Table \ref{tab:errors patch 1} presents the analytical results for the stresses and compares them to the minimum and maximum computed value at every dof of the mesh. The maximum relative error is also given.
\begin{table}[!htp]
\begin{center} \footnotesize
   \begin{tabular}{ | c | c | c | c | c | c | c |}
     \hline 
      stress & $\sigma_{xx}$ & $\sigma_{yy}$ & $\sigma_{xy}$ & $\sigma_{yx}$ & $\mu_x$ & $\mu_y$  \\ \hline
     analytical & $4$ & $4$ & $1.5$ & $1.5$ & $0$ & $0$ \\ \hline
     min computed & $4.00$ &  $4.00$ & $1.50$ & $1.50$ & ME & ME \\ \hline
     max computed & $4.00$ &  $4.00$ & $1.50$ & $1.50$ & ME & ME \\ \hline
     max relative error & $2.63\cdot 10^{-11}\%$ & $6.21\cdot 10^{-11}\%$ & $4.15\cdot 10^{-11}\%$ & $1.04\cdot 10^{-10}\%$ & &  \\ \hline
   \end{tabular}
   \caption{First patch test: Analytical solution, computed stresses and relative error.}
   \label{tab:errors patch 1}
\end{center}
\end{table}
ME means machine error and the relative error on the momenta is not given as the expected value is zero.

\subsubsection{Second patch test}
The solution to the second test is
\begin{equation}
\left\{
\begin{aligned}
& u_x(x,y) = \frac1{G}(x + \frac{y}2), \\
& u_y(x,y) = \frac1{G}(x + y), \\
&\varphi(x,y) = -\frac1{4G}.
\end{aligned}
\right.
\end{equation}
Dirichlet boundary conditions are imposed on the entire boundary.
The volumic load is $f_x = f_y = 0$ and $\mathfrak{c} = 1$.
The condition number with the penalty term \eqref{eq:innner penalty} is $31$ and $28$ without it.
Table \ref{tab:errors patch 2} presents the analytical and computed results.
\begin{table}[!htp]
\begin{center} \footnotesize
   \begin{tabular}{ | c | c | c | c | c | c | c |}
     \hline 
      stress & $\sigma_{xx}$ & $\sigma_{yy}$ & $\sigma_{xy}$ & $\sigma_{yx}$ & $\mu_x$ & $\mu_y$  \\ \hline
     analytical & $4$ & $4$ & $1$ & $2$ & $0$ & $0$ \\ \hline
     min computed & $3.99$ &  $3.98$ & $0.995$ & $1.96$ & $-4.58\cdot 10^{-2}$ & $-5.36\cdot 10^{-2}$ \\ \hline
     max computed & $4.01$ &  $4.01$ & $1.03$ & $2.01$ & $4.23\cdot 10^{-2}$ & $5.42\cdot 10^{-2}$ \\ \hline
     max relative error & $1.58\%$ & $1.53\%$ & $3.51\%$ & $2.35\%$ & &  \\ \hline
   \end{tabular}
   \caption{Second patch test: Analytical solution, computed stresses and relative error.}
   \label{tab:errors patch 2}
\end{center}
\end{table}

\subsubsection{Third patch test}
The solution to the third test is
\begin{equation}
\left\{
\begin{aligned}
& u_x(x,y) = \frac1{G}(x + \frac{y}2), \\
& u_y(x,y) = \frac1{G}(x + y), \\
&\varphi(x,y) = \frac1{G}(\frac14 - x+y)).
\end{aligned}
\right.
\end{equation}
Dirichlet boundary conditions are imposed on the entire boundary.
The volumic load is $f_x = f_y = 1$ and $\mathfrak{c} = 2(x-y)$.
The condition number with the penalty term \eqref{eq:innner penalty} is $112$ and $124$ without it.
Table \ref{tab:errors patch 3} presents the analytical and computed results.
\begin{table}[!htp]
\begin{center} \footnotesize
   \begin{tabular}{ | c | c | c | c | c | c | c |}
     \hline 
      stress & $\sigma_{xx}$ & $\sigma_{yy}$ & $\sigma_{xy}$ & $\sigma_{yx}$ & $\mu_x$ & $\mu_y$  \\ \hline
     analytical & $4$ & $4$ & $1.5-x+y$ & $1.5+x-y$ & $-4\ell^2$ & $4\ell^2$ \\ \hline
     min computed & $3.94$ &  $3.94$ &  &  & $-4.80\cdot 10^{-2}$ & $3.63\cdot 10^{-2}$ \\ \hline
     max computed & $4.06$ &  $4.06$ &  & & $-3.75\cdot 10^{-2}$ & $4.73\cdot 10^{-2}$ \\ \hline
     max relative error & $1.58\%$ & $1.53\%$ & $3.51\%$ & $2.35\%$ & $6.22\%$ & $9.29\%$  \\ \hline
   \end{tabular}
   \caption{Third patch test: Analytical and computed stresses and relative error.}
   \label{tab:errors patch 3}
\end{center}
\end{table}
The minimum and maximum computed values for $\sigma_{xy}$ and $\sigma_{yx}$ are not given for this test because they are irrelevant, as they vary with the position in the domain.

\subsubsection{Results}
Results for finer meshes are not shown here, but the errors decrease with mesh refinement.
The results provided above validate the correct imposition of Dirichlet boundary conditions with the proposed method.
Also, one can notice that the condition number of the rigidity matrix remains of the same order when adding the inner penalty term \eqref{eq:innner penalty}.

\subsection{Numerical results}
\label{sec:numerical results}

\subsubsection{Square plate with a hole}
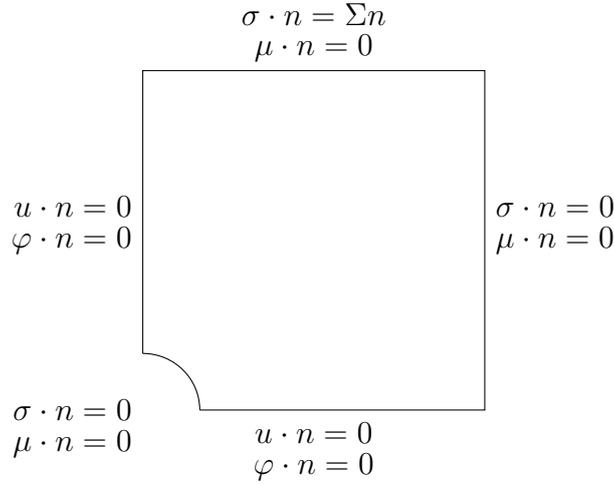
\begin{figure} [!htp]
\begin{center}
\begin{tikzpicture} [scale = 1.5]
\pgfmathsetmacro{\L}{3}
\pgfmathsetmacro{\R}{0.5}

\draw (\R,0) -- (\L,0) -- (\L,\L) -- (0,\L) -- (0,\R);
\draw (\R,0) arc (0:90:\R);

\draw[left] (0,0) node{$\sigma \cdot n = 0$};
\draw[left] (0,-0.3) node{$\mu \cdot n = 0$};
\draw (\L/2,-0.2) node{$u \cdot n = 0$};
\draw (\L/2,-0.5) node{$\varphi \cdot n = 0$};
\draw[left] (0,\L/2+0.3) node{$u \cdot n = 0$};
\draw[left] (0,\L/2) node{$\varphi \cdot n = 0$};
\draw (\L/2,\L+0.5) node{$\sigma \cdot n = \Sigma n$};
\draw (\L/2,\L+0.2) node{$\mu \cdot n = 0$};
\draw[right] (\L,\L/2+0.3) node{$\sigma \cdot n = 0$};
\draw[right] (\L,\L/2) node{$\mu \cdot n = 0$};

\end{tikzpicture}
\caption{Square plate with a hole: problem setup.}
\label{fig:plate with a hole}
\end{center}
\end{figure}

\begin{figure} [!htp]
\subfloat{
\includegraphics[trim=350 50 350 50, clip,width=0.5\textwidth]{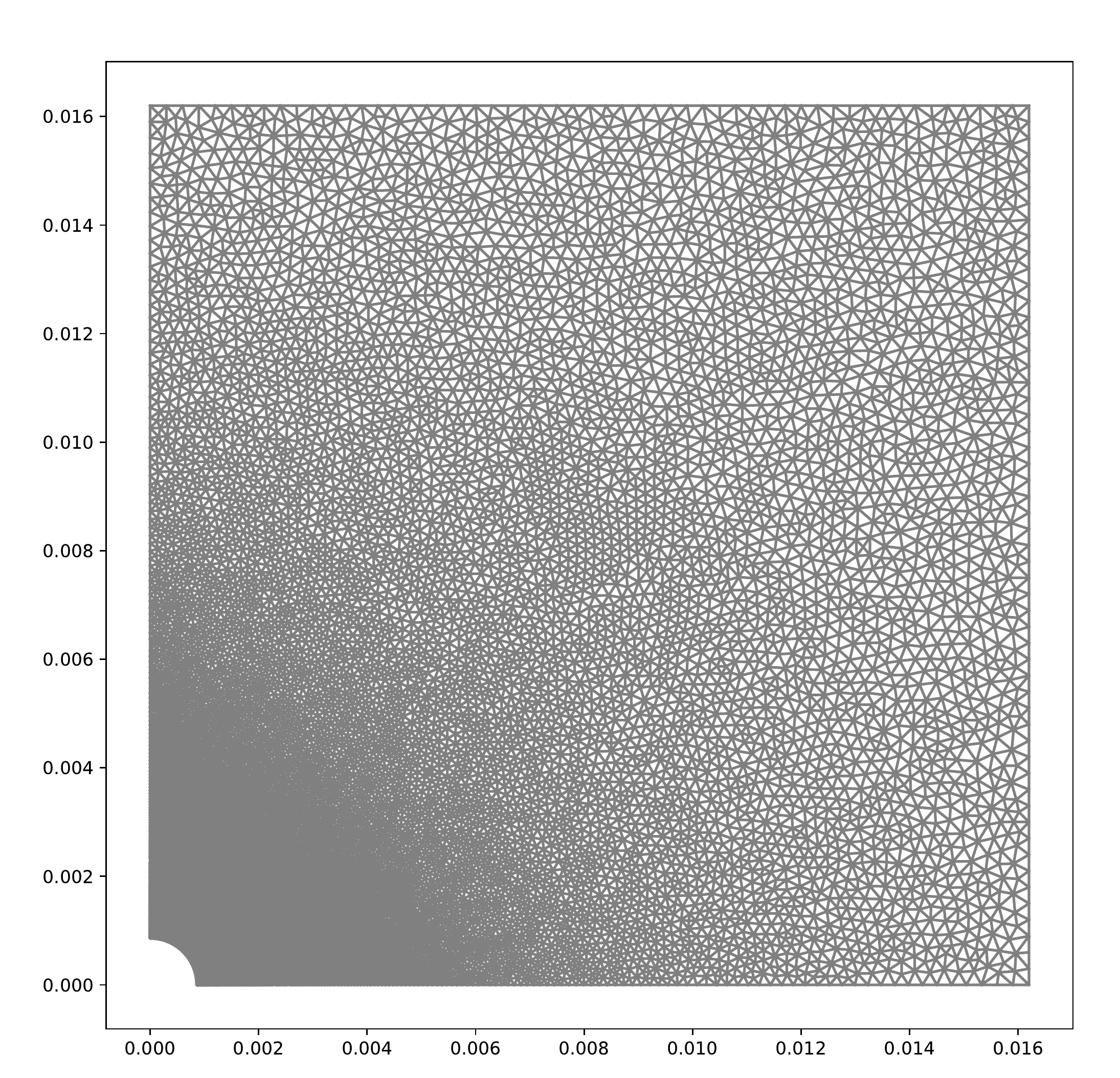}
}
\subfloat{
\includegraphics[trim=350 50 350 50, clip, width=0.5\textwidth]{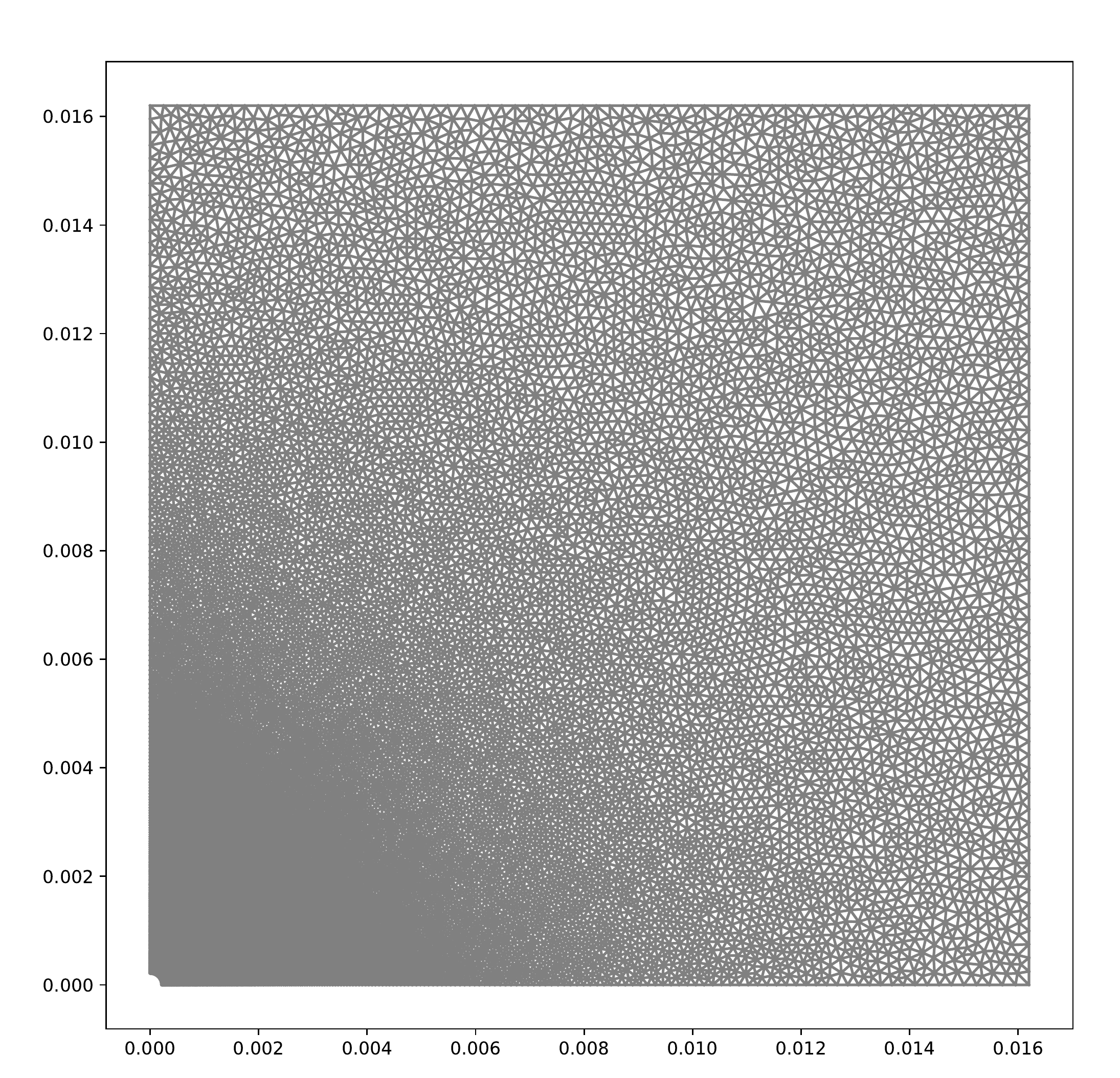}
}
\caption{Plate with a hole: left: mesh for first and second test $r=0.216$mm. right: mesh for third test $r=0.864$mm.}
\label{fig:mesh plate}
\end{figure}

This test case has been inspired by \cite{providas2002finite}.
The domain is a square plate of length $32.4$mm with a hole of radius $r$ in its center. For symmetry reasons, only a quarter of the plate is meshed. A traction with $\Sigma = 1N/m^2$ is imposed on the top surface. Figure \ref{fig:plate with a hole} shows the boundary conditions.
The material parameters are as follows:
$G=10^3$Pa and $\nu = 0.3$. $\ell$ and $a$ as well as the radius $r$ take several values as presented in the following. The meshes in Figure \ref{fig:mesh plate} are used for the test cases.
The mesh for the first two tests is made of $225,816$ dofs and the mesh for the third test is made of $210,867$ dofs.
For the first test case, $r=0.216$mm and $\frac{r}{\ell} = 1.063$ and $a$ varies as in Table \ref{tab:errors 2d 1} which gives the maximal stress at the boundary of the hole as well as the error.
\begin{table}[!htp]
\begin{center}
   \begin{tabular}{ | c | c | c | c | c |}
     \hline 
     $a$ & analytical & computed & error (\%) \\ \hline
     $0$ & $3.000$ &  $2.998$ & $0.1\%$ \\ \hline
     $0.0667$ & $2.849$ & $2.848$ & $0.0\%$ \\ \hline
     $0.3333$ & $2.555$ & $2.555$ & $0.0\%$ \\ \hline
     $1.2857$ & $2.287$ & $2.287$ & $0.0\%$ \\ \hline
     $4.2632$ & $2.158$ & $2.157$ & $0.0\%$ \\ \hline
   \end{tabular}
   \caption{Square plate with a hole: first test, parameter $a$, analytical maximal stress, computed maximal stress and error.}
   \label{tab:errors 2d 1}
\end{center}
\end{table}
For the second test case, $r=0.216$mm and $\frac{r}{\ell} = 10.63$ and $a$ varies as in Table \ref{tab:errors 2d 2} which gives the maximal stress at the boundary of the hole as well as the error.
\begin{table}[!htp]
\begin{center}
   \begin{tabular}{ | c | c | c | c | c |}
     \hline 
     $a$ & analytical & computed & error (\%) \\ \hline
     $0$ & $3.000$ &  $2.998$ & $0.1\%$ \\ \hline
     $0.0667$ & $2.956$ & $2.955$ & $0.0\%$ \\ \hline
     $0.3333$ & $2.935$ & $2.936$ & $0.0\%$ \\ \hline
     $1.2857$ & $2.927$ & $2.929$ & $0.1\%$ \\ \hline
     $4.2632$ & $2.923$ & $2.925$ & $0.1\%$ \\ \hline
   \end{tabular}
   \caption{Square plate with a hole: second test, parameter $a$, analytical maximal stress, computed maximal stress and error.}
   \label{tab:errors 2d 2}
\end{center}
\end{table}
For the third test case, $r=0.864$mm and $a=0.3333$mm and $\ell$ varies as in Table \ref{tab:errors 2d 3} which gives the maximal stress at the boundary of the hole as well as the error.
\begin{table}[!htp]
\begin{center}
   \begin{tabular}{ | c | c | c | c | c |}
     \hline 
     $r/\ell$ & analytical & computed & error (\%) \\ \hline
     $1.0$ & $2.549$ &  $2.566$ & $0.7\%$ \\ \hline
     $2.0$ & $2.641$ &  $2.660$ & $0.7\%$ \\ \hline
     $3.0$ & $2.719$ & $2.740$ & $0.8\%$ \\ \hline
     $4.0$ & $2.779$ & $2.801$ & $0.8\%$ \\ \hline
     $6.0$ & $2.857$ & $2.881$ & $0.8\%$ \\ \hline
     $8.0$ & $2.902$ & $2.927$ & $0.9\%$ \\ \hline
     $10.0$ & $2.929$ & $2.955$ & $0.9\%$ \\ \hline
   \end{tabular}
   \caption{Square plate with a hole: third test, parameter $\ell$, analytical maximal stress, computed maximal stress and error.}
   \label{tab:errors 2d 3}
\end{center}
\end{table}
The results from the three tests show that the proposed method is able to reproduce stress concentration with a satisfactory accuracy.

\subsubsection{Boundary layer effect}
\begin{figure} [!htp]
\centering
\subfloat{
\includegraphics[width=0.5\textwidth]{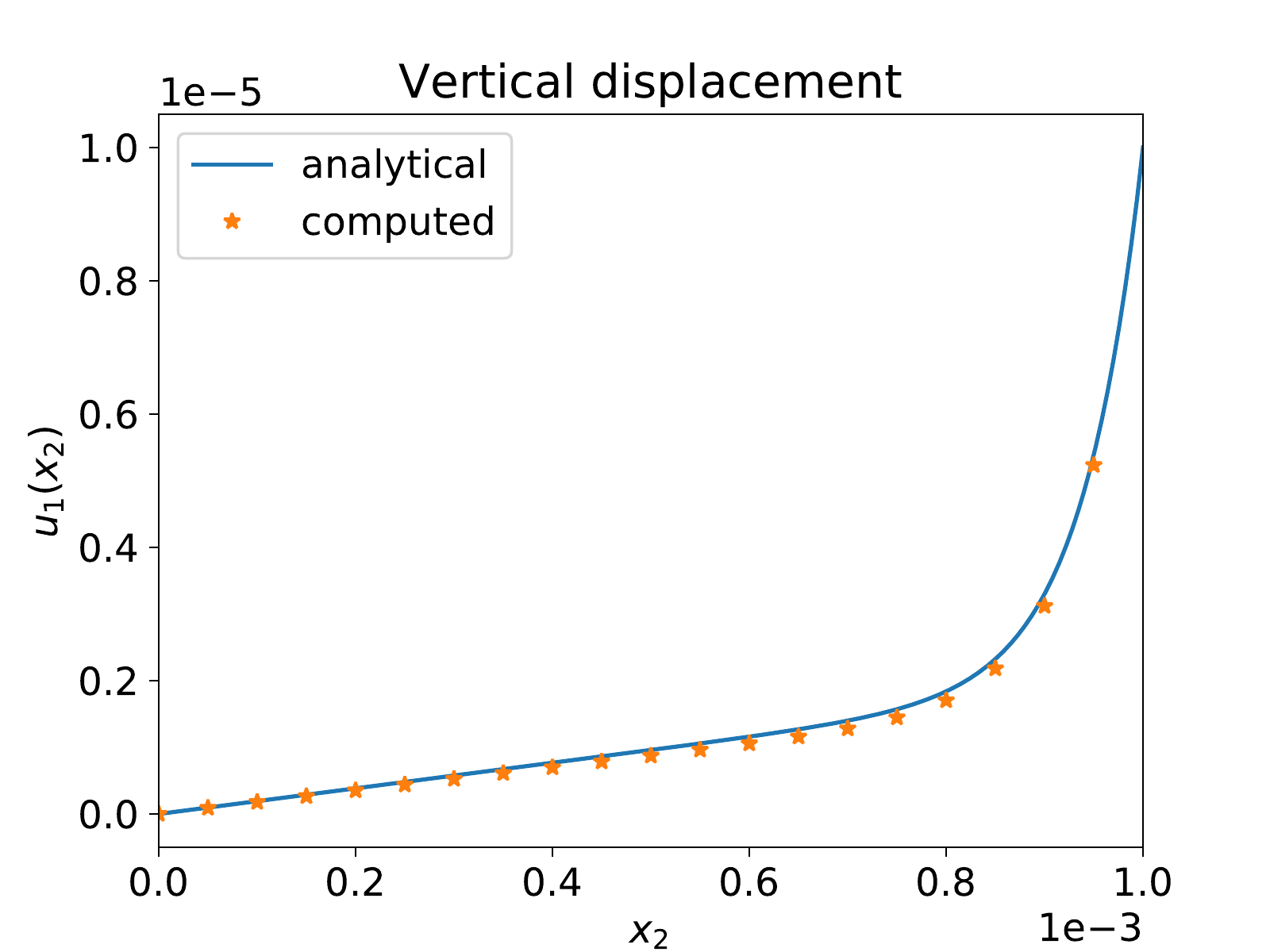}
}
\subfloat{
\includegraphics[width=0.5\textwidth]{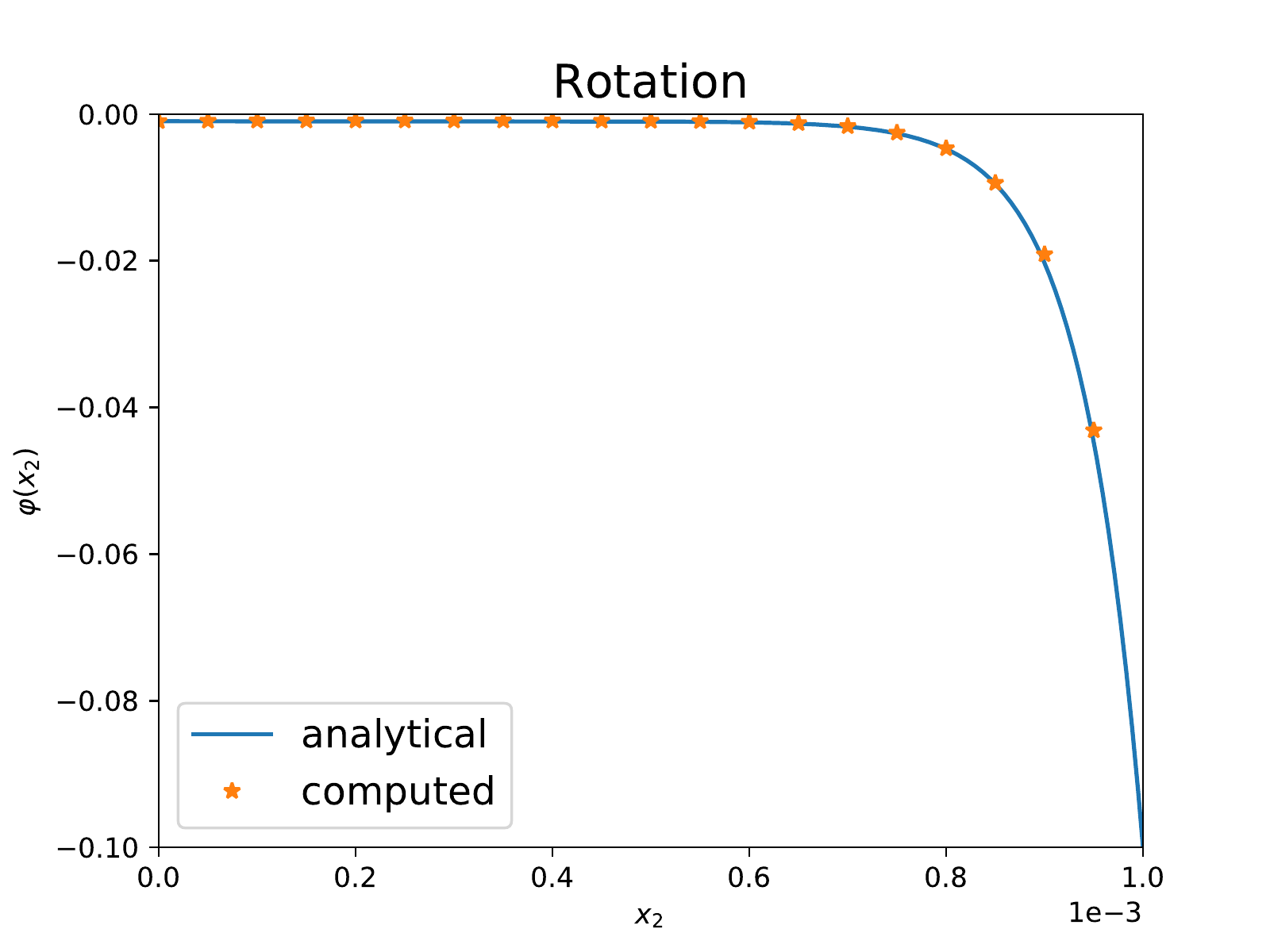}
}
\caption{Boundary layer effect: left: displacement in the $x_1$ direction. right: rotation.}
\label{fig:boundary layer}
\end{figure}

This test case is found in \cite[p. 342]{sulem1995bifurcation}, for the analytical solution and in \cite{rattez2018numerical}, for a numerical implementation.
The domain is a square of size $h=1$mm. The material parameters are $\nu=0$, $G=10$GPa, $a=2$ and $\ell = 5\cdot 10^{-2}$mm using the material laws \eqref{eq:2d material law} and \eqref{eq:2d material law 2}. The boundary conditions are $u_1=0$ and $\varphi=0$ on the bottom surface, $u_1=-0.1$m and $\varphi=0.01\cdot h$ on the top surface. Mirror boundary conditions are imposed on the left and right boundaries. The mesh is a structured triangular mesh made of 50 elements
along the $x_2$-direction and 10 along the $x_1$-direction thus leading to $6,000$ dofs.
Figure \ref{fig:boundary layer} shows the computed values of $u_1$ and $\varphi$ depending on the $x_2$ coordinate compared to the analytical solution available in \cite{rattez2018numerical}.
A similar computation is performed with mixed Lagrange $P^2$-$P^1$ FE on a mesh containing $6,003$ dofs. The maximum relative error on the computed values ploted in Figure \ref{fig:boundary layer} is $9\%$ for the rotation and the vertical displacement for the DEM. It is of $1\%$ for the vertical displacement and the rotation with the FE computation.

\section{Fully discrete scheme}
\label{sec:dynamics}
A Crank--Nicholson time-integration \cite{belytschko1983computational} is use to integrate in time the system of ODEs \eqref{eq:discrete dynamics}.
\subsection{Space-time discrete system}
The time step is chosen as $\Delta t = \frac{T}{2000}$.
The time-interval $(0,T)$ is discretized in $\{0=t_0, \dots , t^n, \dots, t^N=T\}$.
For all $n=1,\ldots,N$, we compute the discrete displacement and rotation $u_h^n$ and $\varphi_h^n$ and the discrete velocity $\dot{u}_h^{n}$ and rotation rate $\dot{\varphi}_h^{n}$ as well as the corresponding accelerations $\ddot{u}_h^{n}$ and $\ddot{\varphi}_h^{n}$.
As, the homogeneous Dirichlet boundary conditions are imposed weakly in the DEM through \eqref{eq:nitsche con} and \eqref{eq:nitsche nsym}, a damping term is added to impose homogeneous Dirichlet boundary conditions on the velocities:
\begin{equation}
\label{eq:damping}
c_h((\dot{u}_h^{n},\dot{\varphi}_h^{n});(v_h,\psi_h))) := \sum_{F \in \mathcal{F}^b_D} \frac{4G}{h_F} \int_F (\dot{u}^{n}_h \cdot v_h + \ell^2 \dot{\varphi}^{n}_h \cdot \psi_h),  \quad \forall (v_h,\psi_h).
\end{equation}
The fully discrete scheme reads as follows: for all $n=1,\ldots,N$, given
$(u_h^n,\varphi_h^n)$, $(\dot{u}_h^{n},\dot{\varphi}_h^{n})$ and $(\ddot{u}_h^{n},\ddot{\varphi}_h^{n})$, compute $(u_h^{n+1},\varphi_h^{n+1})$, $(\dot{u}_h^{n+1},\dot{\varphi}_h^{n+1})$ and $(\ddot{u}_h^{n+1},\ddot{\varphi}_h^{n+1})$ such that
\begin{equation} 
\label{eq:Crank}
\left\{ \begin{aligned}
&u_h^{n+1} = u_h^n + \Delta t \dot{u}_h^{n} + \frac{\Delta t^2}{2} \frac{\ddot{u}_h^{n}+\ddot{u}_h^{n+1}}{2}, \quad \varphi_h^{n+1} = \varphi_h^n + \Delta t \dot{\varphi}_h^{n} + \frac{\Delta t^2}{2} \frac{\ddot{\varphi}_h^{n}+\ddot{\varphi}_h^{n+1}}{2}, \\
&\dot{u}_h^{n+1} = \dot{u}_h^n + \Delta t \frac{\ddot{u}_h^{n}+\ddot{u}_h^{n+1}}{2}, \quad \dot{\varphi}_h^{n+1} = \dot{\varphi}_h^n + \Delta t \frac{\ddot{\varphi}_h^{n}+\ddot{\varphi}_h^{n+1}}{2}, \\
&\ddot{u}_h^{n+1} = \frac4{\Delta t^2} (u_h^{n+1} - u_h^n - \Delta t \dot{u}_h^{n}) - \ddot{u}_h^{n}, \quad \ddot{\varphi}_h^{n+1} = \frac4{\Delta t^2} (\varphi_h^{n+1} - \varphi_h^n - \Delta t \dot{\varphi}_h^{n}) - \ddot{\varphi}_h^{n}, \\
&m_h((\ddot{u}_h^{n+1},\ddot{\varphi}_h^{n+1});(v_h,\psi_h)) + c_h((\dot{u}_h^{n},\dot{\varphi}_h^{n});(v_h,\psi_h))) \\ +
& a_h((u_h^{n+1},\phi_h^{n+1});(v_h,\psi_h))) = L_h(t^{n+1};(v_h,\psi_h)) , \quad\forall (v_h,\psi_h), \\
\end{aligned} \right.
\end{equation}
where $m_h$ is defined in \eqref{eq:mass bilinear form} and $L_h$ is the discrete load linear form (the right-hand side of \eqref{eq:discrete dynamics}).
The initial displacement and rotation $(u_h^0,\varphi_h^0)$ and the initial velocity and rotation rate $(\dot{u}_h^{0}, \dot{\varphi}_h^0)$ are evaluated by
using the values of the prescribed initial displacements and velocities at the cell barycentres. 
Note that, with respect to the DEM developed in \cite{LM_CM_2012}, there is no need to solve a nonlinear problem to compute the rotation at each time-step, which greatly improves numerical efficiency.

\subsection{Numerical results}
\subsubsection{Beam in dynamic flexion}
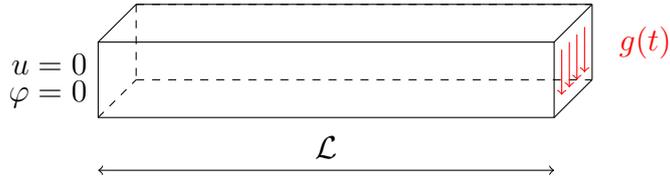
\begin{figure} [!htp]
\begin{center}
\begin{tikzpicture}
\coordinate (a) at (-3,-0.5);
\coordinate (b) at (-3,0.5);
\coordinate (c) at (3,0.5);
\coordinate (d) at (3,-0.5);
\coordinate (ap) at (-2.5,0);
\coordinate (bp) at (-2.5,1.);
\coordinate (cp) at (3.5,1.);
\coordinate (dp) at (3.5,0);

\draw (a) -- (b) -- (c) -- (d) -- cycle;
\draw[dashed] (a) -- (ap);
\draw (b) -- (bp);
\draw[dashed] (bp) -- (ap);
\draw (c) -- (cp);
\draw (d) -- (dp);
\draw (cp) -- (dp);
\draw (bp) -- (cp);
\draw[dashed] (bp) -- (cp);
\draw[dashed] (ap) -- (dp);
\draw [<->] (-3,-1.2) -- (3,-1.2);
\draw (0,-0.9) node{$\mathcal{L}$}; 
\draw (4.2,0.5) node[red] {$g(t)$};
\draw [<-,red] (3.4, 0.1) -- (3.4, 0.7);
\draw [<-,red] (3.3, -0.) -- (3.3, 0.6);
\draw [<-,red] (3.2, -0.1) -- (3.2, 0.5);
\draw [<-,red] (3.1, -0.2) -- (3.1, 0.4);
\draw(-3,0.2) node[left] {$u=0$};
\draw(-3,-0.2) node[left] {$\varphi=0$};
\end{tikzpicture}
\caption{Beam in dynamic flexion: problem setup.}
\label{fig:dynamical test case sketch}
\end{center}
\end{figure}

\begin{figure} [!htp]
\centering
\includegraphics[width=0.5\textwidth]{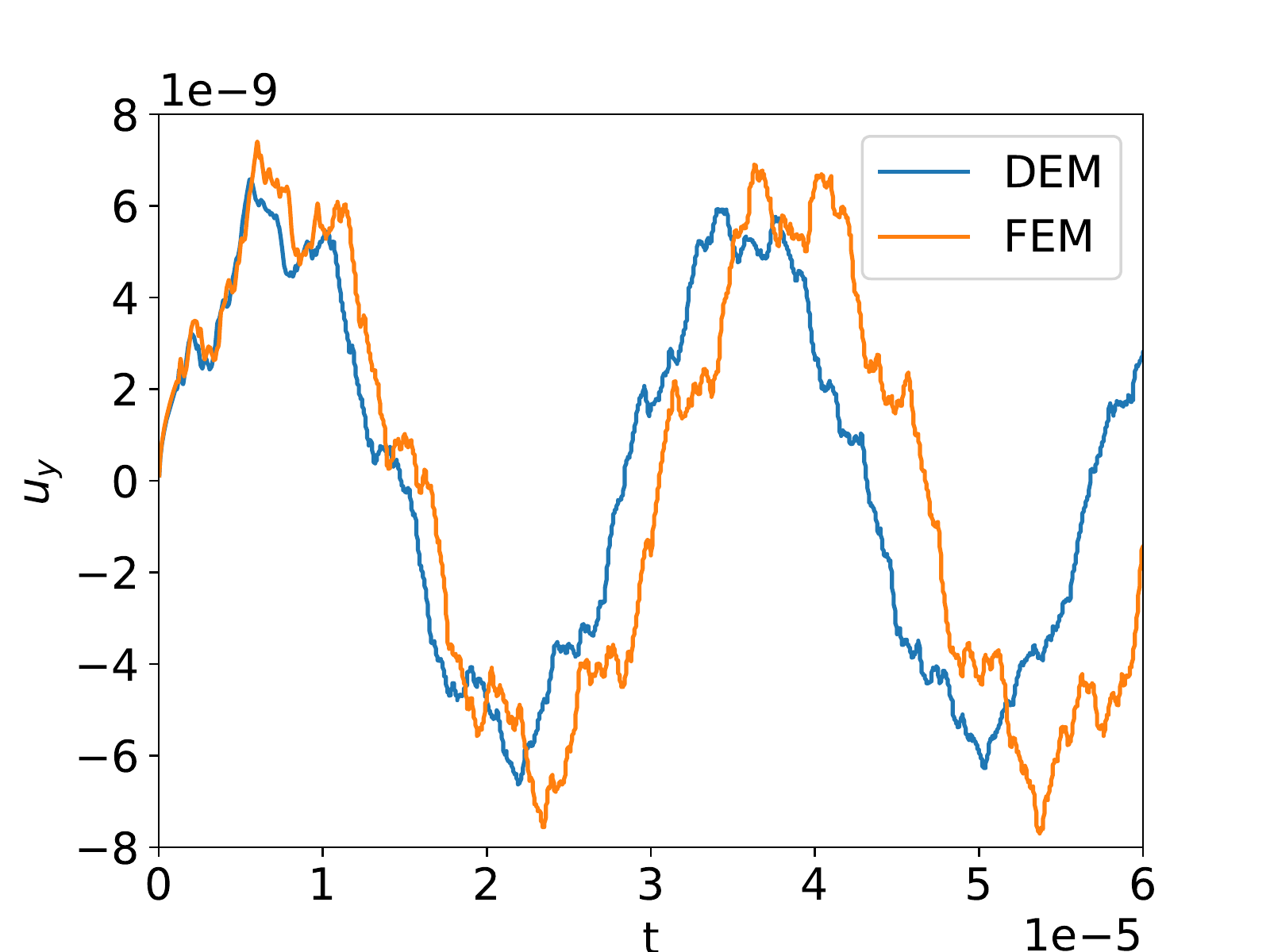}
\caption{Beam in dynamic flexion: displacement of the tip of the beam over simulation time.}
\label{fig:dynamic tip disp}
\end{figure}

This test case has been inspired by a similar from \cite{bleyer2018numericaltours}.
This test case consists of computing the oscillations of a beam of length $\mathcal{L}=1$mm with a square section of $0.04 \times 0.04$mm$^2$. The simulation time is $T=6.3 \cdot 10^{-5}$s. The beam is clamped at one end, it is loaded by a uniform vertical traction $\sigma\cdot n = g(t)$ at the other end, and the four remaining lateral faces are stress free ($\sigma\cdot n=0$ and $\mu \cdot n = 0$). The load term $g(t)$ is defined as
\begin{equation}
g(t) :=
\begin{cases}
-\frac{tE\cdot10^{-6}}{T_c}e_x& \text{for $0 \leq t \leq T_c$}, \\
0&\text{for $T_c \leq t \leq T$},
\end{cases}
\end{equation}
where $T_c = 3.2 \cdot 10^{-8}$s.
Figure \ref{fig:dynamical test case sketch} displays the problem setup.
The material parameters have been taken from \cite{RATTEZ201854}. The bulk modulus is $K = 16.67$GPa and the shear moduli are $G = 10$GPa and $G_c = 5$GPa.
The characteristic size of the microstructure is taken as $\ell = \frac{\mathcal{L}}{100}$. We also take $L = G\ell^2$ and $M = M_c = \frac52 Gl^2$.
The density is $\rho = 2500\mathrm{kg {\cdot} m^{-3}}$ and the inertia is taken as $I = \frac25 \ell^2$ following \cite{RATTEZ201854}, with the assumption that the micro-structure of the material is composed of balls.
The reference solution is a $P^2$-$P^1$ Lagrange FEM coupled to a Crank--Nicholson time-integration \cite{belytschko1983computational}. The DEM is integrated according to Equation \eqref{eq:Crank}.
The DEM computation contains $23,040$ dofs and the FEM computation $3,642$ dofs.
Figure \ref{fig:dynamic tip disp} shows the displacement of the beam tip $u_h(\mathcal{L}, 0.05, 0.) \cdot e_y$ over $[0,T]$ for the two computations.
As expected, the two methods deliver similar results.

\subsubsection{Lamb's problem}
\begin{figure} [!htp]
\subfloat{
\includegraphics[width=0.8\textwidth]{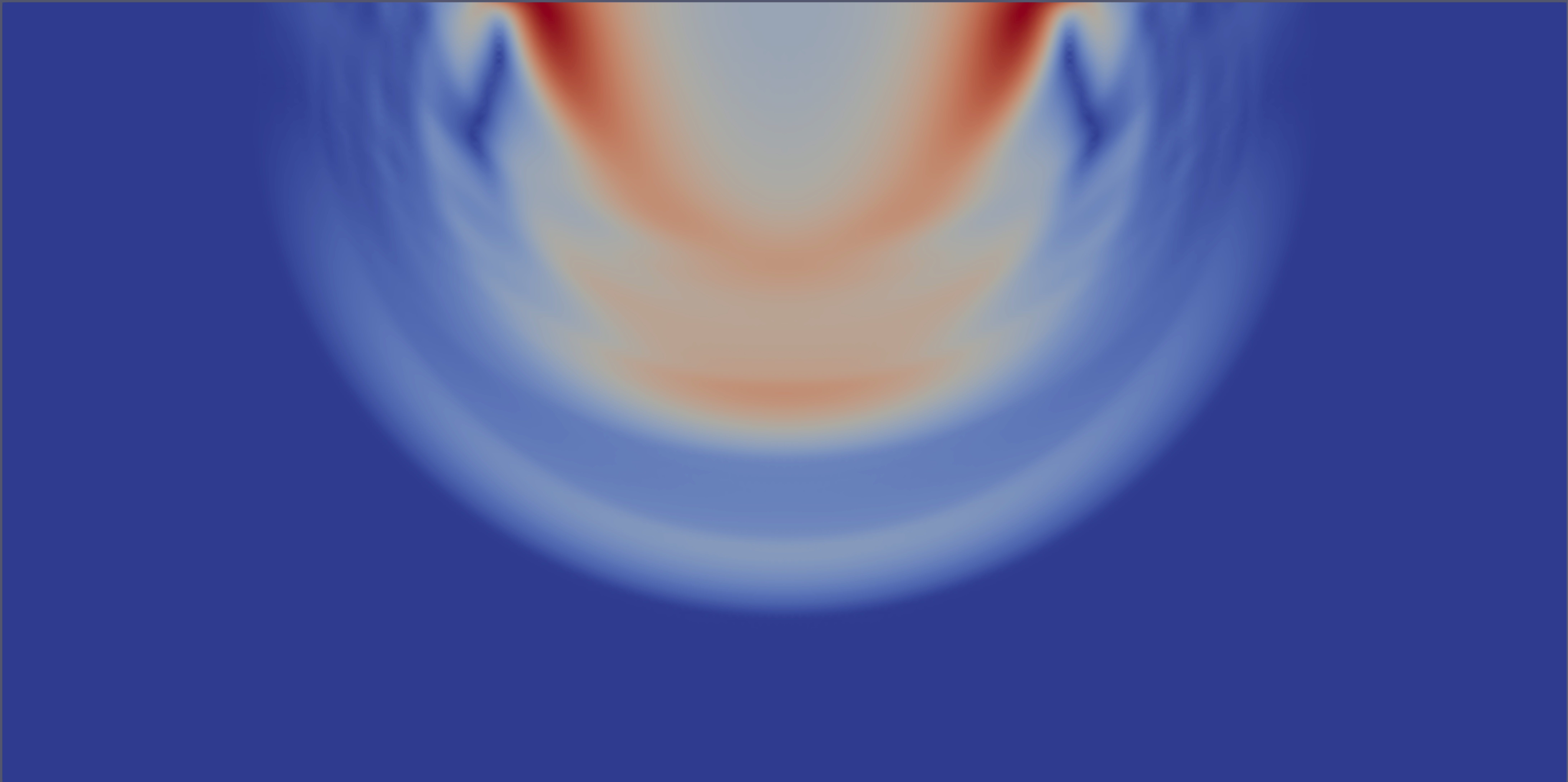}
}
\subfloat{
\includegraphics[width=0.065\textwidth]{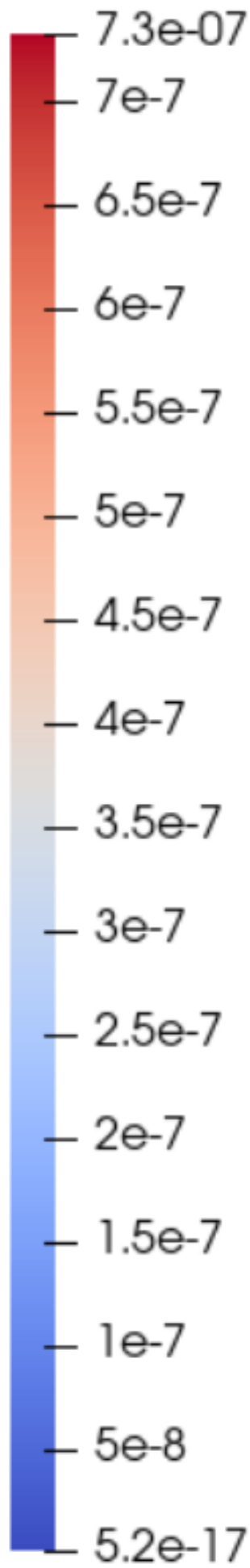}
}\\
\subfloat{
\includegraphics[width=0.8\textwidth]{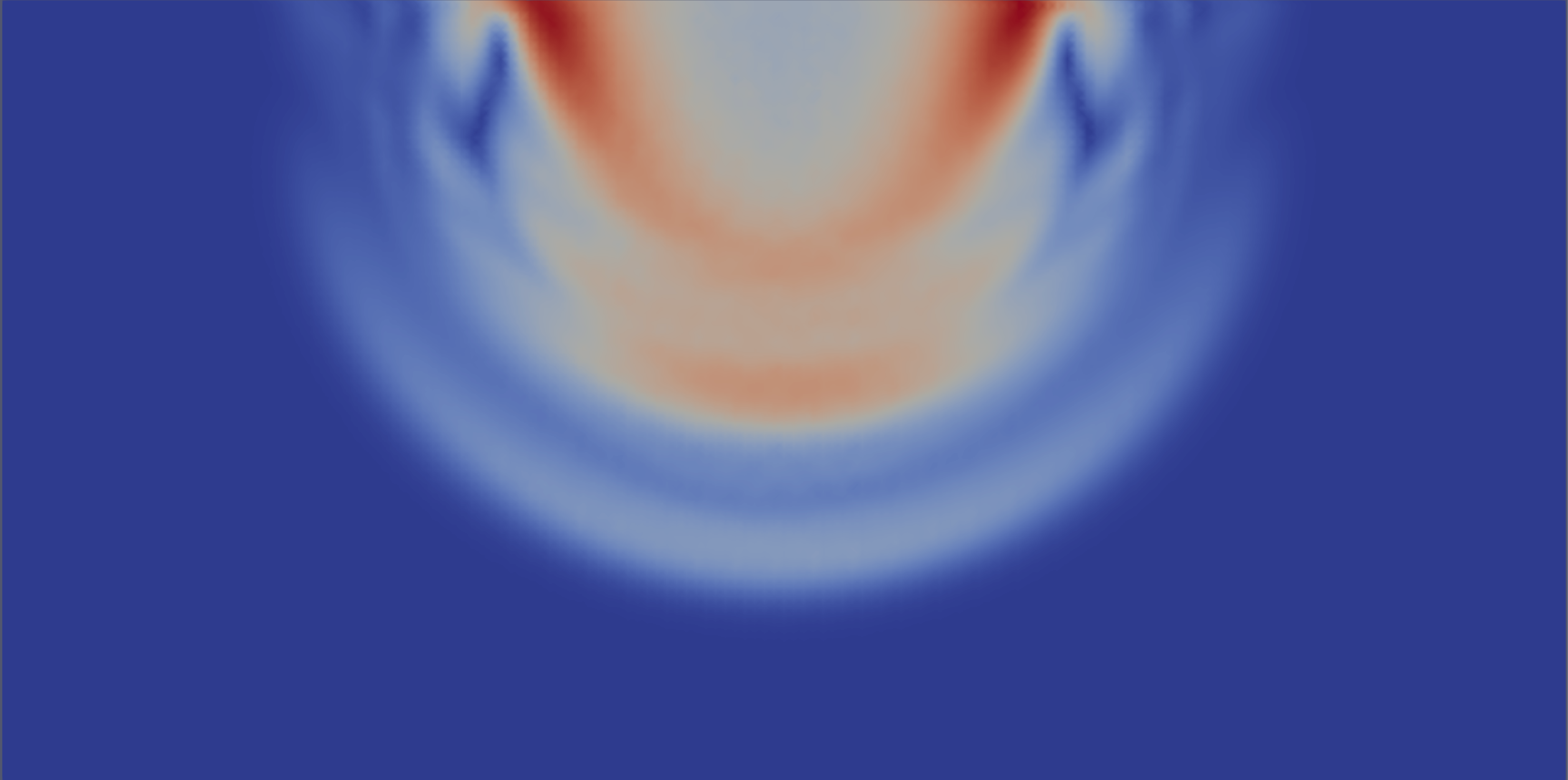}
}
\subfloat{
\includegraphics[width=0.065\textwidth]{Lamb/colorbar.pdf}
}
\caption{Lamb's problem: velocity magnitude at $t=0.2$s. top: FEM, bottom: DEM.}
\label{fig:Lamb}
\end{figure}

Lamb's problem \cite{lamb1904propagation} is a classical test case used to assert the capacity of a numerical method is reproduce the propagation of seismic waves.
Following \cite{mariotti2007lamb}, we consider a rectangular domain of size $2 \times 1$km$^2$.
A source modelled by a Ricker pulse of central frequency $14.5$Hz is placed $100$m below the top surface in the middle of the rectangle. Homogeneous Neumann boundary conditions are imposed on the entire boundary.
Following \cite{LM_CM_2012}, the material parameters are taken as $G = G_c = 7.52$GPa, $\lambda = 3.76$GPa and $\ell = \frac{h}{\sqrt{2}}$, where $h$ is the size of the mesh, and the inertia is taken as $I = \frac{\ell^2}{6}$. 
Such a choice for $\ell$ can look baffling if one considers that the computed material is indeed a Cosserat continuum. However, if one considers a variational DEM approach to compute seismic waves, then such a choice is backed by the litterature \cite{LM_CM_2012}.
The DEM computation is performed with $153,600$ dofs and a time-step $\Delta t = 5.0 \cdot 10^{-5}$s. The reference solution is a P$^2$-P$^1$ Lagrange FEM with $361,503$ dofs and a similar time-step.
Figure \ref{fig:Lamb} shows the magnitude of the velocity vector at $t=0.2$s for the DEM and the FEM computation.
One can observe the propagation of three waves. First a compression P-wave and then a shear S-wave propagate inside the domain. Finally, a Rayleigh wave propagates on the upper boundary (in red).
The results given by the two methods coincide strongly and confirm the pertinence of using a Cosserat material law and DEM to compute seismic waves as proposed in \cite{mariotti2007lamb,LM_CM_2012}.

\section{Conclusion}
\label{sec:conclusion}
In this article, a variational DEM has been introduced which features only cell unknowns for the displacement and the micro-rotation. A cellwise gradient reconstruction is used to obtain cellwise constant strains and stresses using the formalism of Cosserat materials.
An interpretation of the method as a DEM is presented in which the forces exerted by every facet (or link) between two cells (or discrete elements) are explicitly given as functions of the cellwise constant reconstructed stresses.
The method is proved to give satisfactory results on many different test cases in both two and three space dimensions and both in statics and dynamics.
Further work could include extending the present formalism to nonlinear material laws in small strains \cite{RATTEZ201854} and then finite strains \cite{neff2006finite}. It could also include computing nonlinear dynamic evolutions \cite{neff2007well}.
Also, the present formalism could be extended to Gyro-continua \cite{brocato2001gyrocontinua} to have a real rotation matrix in each cell.

\section*{Acknowledgements}
The author would like to thank A. Ern from Inria and Ecole Nationale des Ponts et Chaussées for stimulating discussions and L. Monasse from Inria for carefully proof-reading this manuscript.
The author would also like to thank the anonymous reviewers for their contributions which helped substantially improve this paper.

\section*{Funding}
Not applicable.

\section*{Conflicts of interest/Competing interests}
The author has no conflict of interest or competing interests.

\section*{Availability of data and material}
Not applicable.

\section*{Code availability}
\url{https://github.com/marazzaf/DEM_cosserat.git}

\bibliographystyle{plain} 
\bibliography{biblio}

\end{document}